\documentclass[12pt]{article}%
\usepackage{amsmath}
\usepackage{amsfonts}
\usepackage{amssymb}
\usepackage{graphicx}%
\providecommand{\U}[1]{\protect\rule{.1in}{.1in}}
\newtheorem{theorem}{Theorem}

\newtheorem{conjecture}[theorem]{Conjecture}
\newtheorem{corollary}[theorem]{Corollary}

\newtheorem{definition}[theorem]{Definition}

\newtheorem{exercise}[theorem]{Exercise}
\newtheorem{lemma}[theorem]{Lemma}

\newtheorem{problem}[theorem]{Problem}
\newtheorem{proposition}[theorem]{Proposition}

\newenvironment{proof}[1][Proof]{\noindent\textbf{#1.} }{\ \rule{0.5em}{0.5em}}
\begin{document}

\title{Working with 2s and 3s }
\author{Diego Dominici \thanks{e-mail: dominicd@newpaltz.edu}\\Department of Mathematics\\State University of New York at New Paltz\\1 Hawk Dr.\\New Paltz, NY 12561-2443\\USA\\Phone: (845) 257-2607\\Fax: (845) 257-3571 }
\maketitle

\begin{abstract}
We establish an equivalent condition to the validity of the Collatz
conjecture, using elementary methods. We derive some conclusions and show
several examples of our results. We also offer a variety of exercises,
problems and conjectures.

\end{abstract}

\section{Introduction}

The \textit{Collatz conjecture} (also known as the $3n+1$ conjecture, Ulam's
conjecture, the Syracuse problem, Kakutani's problem, Hasse's algorithm, etc.)
was first proposed by Lothar Collatz in 1937 \cite{Collatz}. In terms of the
function $T(n),$ defined by%
\begin{equation}
T(n)=\left\{
\begin{array}
[c]{c}%
\frac{3n+1}{2},\quad n\equiv1(\operatorname{mod}\ 2)\\
\frac{n}{2},\quad n\equiv0(\operatorname{mod}\ 2)
\end{array}
\right.  ,\quad n\in%
\mathbb{N}
, \label{T}%
\end{equation}
the conjecture claims that for all natural numbers $n,$ there exists a natural
number $k$ such that
\[
T^{\left(  k\right)  }\left(  n\right)  =\underset{k\text{ times}}%
{\underbrace{T\circ T\circ\cdots\circ T}}(n)=1.
\]
For example, we have%
\begin{align}
T(3)  &  =5,\ T^{\left(  2\right)  }(3)=8,\ T^{\left(  3\right)
}(3)=4,\ T^{\left(  4\right)  }(3)=2,\ T^{\left(  5\right)  }(3)=1,\nonumber\\
T(7)  &  =11,\ T^{\left(  2\right)  }(7)=17,\ T^{\left(  3\right)
}(7)=26,\ T^{\left(  4\right)  }(7)=13,\ T^{\left(  5\right)  }%
(7)=20,\ \label{ExT}\\
T^{\left(  6\right)  }(7)  &  =10,\ T^{\left(  7\right)  }(7)=5,\ T^{\left(
8\right)  }(7)=8,\ T^{\left(  9\right)  }(7)=4,\ T^{\left(  10\right)
}(7)=2,\ T^{\left(  11\right)  }(7)=1.\nonumber
\end{align}
We define $T^{\left(  \infty\right)  }\left(  n\right)  =1.$

\begin{exercise}
Prove that if $\forall n\in%
\mathbb{N}
$ $\exists$ $k\in%
\mathbb{N}
$ such that $T^{\left(  k\right)  }\left(  n\right)  <n,$ then the Collatz
conjecture is true. The number $k$ is called the \textit{stopping time} of
$n.$
\end{exercise}

As of February 2007, the Collatz conjecture has been verified for numbers up
to $13\times2^{58}=\allowbreak3,746\,,994,\,889,\,\allowbreak972,\,252,\,672$
\cite{Oliveira}. However, the general case remains open.

Introducing the \textit{total stopping time }function $\sigma_{\infty}(n),$
defined by $\sigma_{\infty}(1)=0$ and%
\[
\sigma_{\infty}(n)=\inf\left\{  k\in%
\mathbb{N}
\cup\left\{  \infty\right\}  \ |\ T^{\left(  k\right)  }\left(  n\right)
=1\right\}  ,\quad n\geq2,
\]
we can reformulate the Collatz conjecture as
\begin{equation}
\mathfrak{C}=%
\mathbb{N}
, \tag{C1}%
\end{equation}
where
\begin{equation}
\mathfrak{C}=\left\{  n\in%
\mathbb{N}
\ |\ \sigma_{\infty}(n)<\infty\right\}  . \label{C}%
\end{equation}
From (\ref{ExT}), we have%
\[%
\begin{tabular}
[c]{|c|c|c|c|c|c|c|c|}\hline
$n$ & $2$ & $3$ & $4$ & $5$ & $6$ & $7$ & $8$\\\hline
$\sigma_{\infty}(n)$ & $1$ & $5$ & $2$ & $4$ & $6$ & $11$ & $3$\\\hline
\end{tabular}
\ \ .
\]

\begin{exercise}
Find $\sigma_{\infty}(n)$ for $9\leq n\leq100.$

Hint: The web page http://www.numbertheory.org/php/collatz.html contains an
implementation which allows the computation of $\sigma_{\infty}(n)$ for large
values of $n.$
\end{exercise}

One could consider the inverse problem and try to characterize the sets
$S_{k},$ defined by $S_{0}=\left\{  1\right\}  $ and%
\begin{equation}
S_{k}=\left\{  n\in%
\mathbb{N}
\ |\ \sigma_{\infty}(n)=k\right\}  ,\quad k\geq1. \label{Sk}%
\end{equation}
The first few $S_{k}$ are
\begin{align}
S_{1}  &  =\left\{  2\right\}  ,\ S_{2}=\left\{  4\right\}  ,\ S_{3}=\left\{
8\right\}  ,\ S_{4}=\left\{  5,16\right\}  ,\ S_{5}=\left\{  3,10,32\right\}
,\label{Sk1}\\
\ S_{6}  &  =\left\{  6,20,21,64\right\}  ,\ S_{7}=\left\{
12,13,40,42,128\right\}  ,\ S_{8}=\left\{  24,26,80,84,85,256\right\}
.\nonumber
\end{align}
It is clear from (\ref{Sk}) that $2^{k}\in S_{k}$ $\forall k\in%
\mathbb{N}
_{0},$ where $%
\mathbb{N}
_{0}=%
\mathbb{N}
\cup\left\{  0\right\}  .$ In terms of the sets $S_{k},$ the Collatz
conjecture reads%
\begin{equation}%
{\displaystyle\bigcup\limits_{k=0}^{\infty}}
S_{k}=%
\mathbb{N}
. \tag{C2}%
\end{equation}

\begin{exercise}
Compute $S_{k}$ for $9\leq k\leq100.$

Hint: Consider the inverse map $T^{-1}:%
\mathbb{N}
\rightarrow\mathbb{P}\left(
\mathbb{N}
\right)  ,$ given by
\[
T^{-1}(n)=\left\{
\begin{array}
[c]{c}%
\left\{  2n\right\}  ,\quad n\equiv0,1(\operatorname{mod}\ 3)\\
\left\{  2n,\frac{1}{3}\left(  2n-1\right)  \right\}  ,\quad n\equiv
2(\operatorname{mod}\ 3)
\end{array}
\right.  .
\]

\end{exercise}

The sequence of natural numbers $\left\{  x_{n}^{(m)},\ n\geq0\right\}  $,
defined by $x_{0}^{(m)}=m$ and%
\begin{equation}
x_{n+1}^{\left(  m\right)  }=T\left(  x_{n}^{(m)}\right)  ,\quad0\leq n,
\label{recu}%
\end{equation}
is called the \emph{trajectory or forward orbit} of $m\in%
\mathbb{N}
$. From (\ref{ExT}), we have%
\begin{align*}
\left\{  x_{n}^{(2)}\right\}   &  =\left\{  2,1\right\}  ,\quad\left\{
x_{n}^{(3)}\right\}  =\left\{  3,5,8,4,2,1\right\}  ,\quad\left\{  x_{n}%
^{(4)}\right\}  =\left\{  4,2,1\right\}  ,\\
\left\{  x_{n}^{(5)}\right\}   &  =\left\{  5,8,4,2,1\right\}  ,\quad\left\{
x_{n}^{(6)}\right\}  =\left\{  6,3,5,8,4,2,1\right\}  ,\\
\quad\left\{  x_{n}^{(7)}\right\}   &  =\left\{
7,11,17,26,13,20,10,5,8,4,2,1\right\}  ,\quad\left\{  x_{n}^{(8)}\right\}
=\left\{  8,4,2,1\right\}  .
\end{align*}

\begin{exercise}
Find $\left\{  x_{n}^{(m)}\right\}  $ for $9\leq m\leq100.$
\end{exercise}

Using the sequences $\left\{  x_{n}^{(m)}\right\}  $ we can restate Collatz's
conjecture as%
\begin{equation}%
{\displaystyle\bigcap\limits_{m=2}^{\infty}}
\left\{  x_{n}^{(m)}\right\}  =\left\{  2,1\right\}  . \tag{C3}%
\end{equation}

We can also consider higher order recurrences, i.e., instead of (\ref{recu}),
use
\[
x_{n+i}^{\left(  m\right)  }=T^{\left(  i\right)  }\left(  x_{n}^{(m)}\right)
,\quad0\leq n,
\]
where
\begin{equation}
T^{\left(  i\right)  }\left(  x\right)  =f_{i,j}(x),\quad\text{if }x\equiv
j\left(  \operatorname{mod}2^{i}\right)  ,\quad0\leq j\leq2^{i}-1.
\label{hrecu}%
\end{equation}
For $i=1,2,3,$ we have%

\begin{align*}
f_{1,0}(x)  &  =\frac{x}{2},\quad f_{1,1}(x)=\frac{3x+1}{2}\\
f_{2,0}(x)  &  =\frac{x}{4},\quad f_{2,1}(x)=\frac{3x+1}{4},\quad
f_{2,2}(x)=\frac{3x+2}{4},\quad f_{2,3}(x)=\frac{9x+5}{4}\\
f_{3,0}(x)  &  =\frac{x}{8},\quad f_{3,1}(x)=\frac{9x+7}{8},\quad
f_{3,2}(x)=\frac{3x+2}{8},\quad f_{3,3}(x)=\frac{9x+5}{8}\\
f_{3,4}(x)  &  =\frac{3x+4}{8},\quad f_{3,5}(x)=\frac{3x+1}{8},\quad
f_{3,6}(x)=\frac{9x+10}{8},\quad f_{3,7}(x)=\frac{27x+19}{8}.
\end{align*}

\begin{exercise}
Prove that if the sequence $\left\{  3k+4\ \right\}  \subset\mathfrak{C},$
then the Collatz conjecture is true.
\end{exercise}

In terms of (\ref{hrecu}), the Collatz conjecture reads%
\begin{equation}
\forall n\in%
\mathbb{N}
\ \exists\ m\in%
\mathbb{N}
\text{ \ such that }n\equiv k\left(  \operatorname{mod}2^{m}\right)  \text{
\ and }\frac{d}{dx}f_{m,k}(x)<1.\tag{C4}%
\end{equation}
For example, we have $11\equiv11\left(  \operatorname{mod}2^{5}\right)  $ and%
\[
f_{5,10}(x)=\frac{3x+2}{32},\quad f_{5,11}(x)=\frac{27x+23}{32}.
\]
Thus,%
\[
x_{5}^{\left(  11\right)  }=f_{5,11}(11)=10,\quad x_{10}^{\left(  11\right)
}=f_{5,10}(10)=1.
\]

The literature on the Collatz conjecture is vast and growing rapidly. Rather
than attempting to cover it, we refer the reader to the excellent survey
papers \cite{Lagarias 1} and \cite{Lagarias 2}.

\section{Representation of natural numbers}

Let the sets $\Lambda_{m}$ be defined by $\Lambda_{m}=\left\{  2^{m}\right\}
,\quad0\leq m\leq3$ and%
\[
\Lambda_{m}=\left\{  \left.  n\in%
\mathbb{N}
\ \right\vert \ n=\frac{2^{m}}{3^{l}}-%
{\displaystyle\sum\limits_{k=1}^{l}}
\frac{2^{b_{k}}}{3^{k}}\right\}  ,\quad m\geq4,
\]
for some $m,l,b_{1},\ldots,b_{l}\in%
\mathbb{N}
_{0},$ with
\[
0\leq l\leq m-3\ \ \text{and }\ 0\leq b_{1}<b_{2}<\cdots<b_{l}\leq m-4.
\]
The first few $\Lambda_{m}$ are%
\begin{align*}
\Lambda_{4}  &  =\left\{  \frac{2^{4}}{3^{0}},\ \frac{2^{4}}{3^{1}}%
-\frac{2^{0}}{3^{1}}\right\}  ,\\
\Lambda_{5}  &  =\left\{  \frac{2^{5}}{3^{0}},\ \frac{2^{5}}{3^{1}}%
-\frac{2^{1}}{3^{1}},\ \frac{2^{5}}{3^{2}}-\left(  \frac{2^{0}}{3^{1}}%
+\frac{2^{1}}{3^{2}}\right)  \right\} \\
\Lambda_{6}  &  =\left\{  \frac{2^{6}}{3^{0}},\ \frac{2^{6}}{3^{1}}%
-\frac{2^{0}}{3^{1}},\ \frac{2^{6}}{3^{1}}-\frac{2^{2}}{3^{1}},\ \frac{2^{6}%
}{3^{2}}-\left(  \frac{2^{1}}{3^{1}}+\frac{2^{2}}{3^{2}}\right)  \right\}  .\
\end{align*}
Using the $\left(  l+2\right)  -$tuple $\left(  l,b_{1},\ldots,b_{l},m\right)
$ to represent the number $\frac{2^{m}}{3^{l}}-%
{\displaystyle\sum\limits_{k=1}^{l}}
\frac{2^{b_{k}}}{3^{k}},$ we can write%
\begin{align}
\Lambda_{4}  &  =\left\{  \left(  0,4\right)  ,\ \left(  1,0,4\right)
\right\}  ,\quad\Lambda_{5}=\left\{  \left(  0,5\right)  ,\ \left(
1,1,5\right)  ,\ (2,0,1,5)\right\}  ,\quad\nonumber\\
\Lambda_{6}  &  =\left\{  \left(  0,6\right)  ,\ \left(  1,0,6\right)
\ ,\left(  1,2,6\right)  ,\ (2,1,2,6)\right\}  ,\label{Lambdak}\\
\Lambda_{7}  &  =\left\{  \left(  0,7\right)  ,\ \left(  1,1,7\right)
\ ,\left(  1,3,7\right)  ,\ (2,0,3,7),\ (2,2,3,7)\right\}  ,\nonumber\\
\Lambda_{8}  &  =\left\{  \left(  0,8\right)  ,\ \left(  1,0,8\right)
,\ \left(  1,2,8\right)  \ ,\left(  1,4,8\right)
,\ (2,1,4,8),\ (2,3,4,8)\right\}  .\ \nonumber
\end{align}

\begin{exercise}
Compute $\Lambda_{m}$ for $9\leq m\leq100.$

Hint: (a) If $\left(  v_{1},v_{2},\ldots,v_{l+2}\right)  \in\Lambda_{m},$ then
$\left(  v_{1},v_{2}+1,\ldots,v_{l+2}+1\right)  \in\Lambda_{m+1}.$

(b) $\left(  1,0,2m\right)  \in\Lambda_{2m}$ for all $m\geq2.$
\end{exercise}

Comparing (\ref{Sk1}) with (\ref{Lambdak}), it seems that $S_{m}=\Lambda_{m}.$
The next results will show this to be true$.$

\begin{lemma}
For all $m\in%
\mathbb{N}
_{0},$ we have%
\begin{equation}
T(\Lambda_{m+1})\subset\Lambda_{m}. \label{TLambda}%
\end{equation}

\end{lemma}

\begin{proof}
Let $n\in\Lambda_{m+1}.$ Then,%
\[
n=\frac{2^{m+1}}{3^{l}}-%
{\displaystyle\sum\limits_{k=1}^{l}}
\frac{2^{b_{k}}}{3^{k}}%
\]
and%
\[
T(n)=\frac{2^{m}}{3^{l}}-%
{\displaystyle\sum\limits_{k=1}^{l}}
\frac{2^{b_{k}-1}}{3^{k}},
\]
if $b_{1}>0$ ($n$ \ even) or%
\begin{align*}
T(n)  &  =\frac{1}{2}\left(  \frac{2^{m+1}}{3^{l-1}}-%
{\displaystyle\sum\limits_{k=1}^{l}}
\frac{2^{b_{k}}}{3^{k-1}}+1\right)  =\frac{1}{2}\left(  \frac{2^{m+1}}%
{3^{l-1}}-%
{\displaystyle\sum\limits_{k=2}^{l}}
\frac{2^{b_{k}}}{3^{k-1}}\right) \\
&  =\frac{1}{2}\left(  \frac{2^{m+1}}{3^{l-1}}-%
{\displaystyle\sum\limits_{k=1}^{l-1}}
\frac{2^{b_{k+1}}}{3^{k}}\right)  =\frac{2^{m}}{3^{l-1}}-%
{\displaystyle\sum\limits_{k=1}^{l-1}}
\frac{2^{b_{k+1}-1}}{3^{k}},
\end{align*}
if $b_{1}=0$ ($n$ \ odd). In either case, $T(n)\in\Lambda_{m}.$
\end{proof}

\begin{lemma}
For all $m\in%
\mathbb{N}
_{0},$ we have%
\begin{equation}
\Lambda_{m}\subset S_{m}. \label{lemma2}%
\end{equation}

\end{lemma}

\begin{proof}
We use induction on $m.$ The case of $m=0$ is clearly true, since $\Lambda
_{0}=\left\{  1\right\}  =S_{0}.$

Assuming (\ref{lemma2}) to be true for $m$, let $n\in\Lambda_{m+1}.$ From
(\ref{TLambda}) we have $T(n)\in\Lambda_{m}$ and therefore $\sigma_{\infty
}\left[  T(n)\right]  =m.$ Thus, $\sigma_{\infty}(n)=m+1$ and the result follows.
\end{proof}

\begin{exercise}
Show that
\begin{equation}
T(n)=\left(  n+\frac{1}{2}\right)  \sin^{2}\left(  \frac{\pi}{2}n\right)
+\frac{n}{2}. \label{T1}%
\end{equation}

\end{exercise}

The other inclusion is also true.

\begin{theorem}
\label{Th}For all $m\in%
\mathbb{N}
_{0},$
\[
S_{m}\subset\Lambda_{m}.
\]

\end{theorem}

\begin{proof}
Clearly, $S_{m}=\left\{  2^{m}\right\}  =\Lambda_{m},\quad0\leq m\leq3.$

Let $m\geq4$ and $s\in S_{m}.$ Using (\ref{T1}) we can write the recurrence
(\ref{recu}) as%
\begin{equation}
x_{n+1}^{\left(  s\right)  }=\left[  x_{n}^{(s)}+\frac{1}{2}\right]
\theta_{n}+\frac{x_{n}^{(s)}}{2},\quad x_{0}^{(s)}=s, \label{recu2}%
\end{equation}
where%
\begin{equation}
\theta_{n}=\sin^{2}\left[  \frac{\pi}{2}x_{n}^{(s)}\right]  ,\text{ \ i.e.,
\ }x_{n}^{(s)}\equiv\theta_{n}^{(s)}\left(  \operatorname{mod}\ 2\right)  .
\label{theta}%
\end{equation}
Assuming $\left\{  \theta_{n}\right\}  $ to be a known sequence, the solution
of (\ref{recu2}) is \cite{Agarwal}%
\[
x_{n}^{(s)}=2^{-n}%
{\displaystyle\prod\limits_{j=0}^{n-1}}
\left(  2\theta_{j}^{(s)}+1\right)  \left(  s+%
{\displaystyle\sum\limits_{k=0}^{n-1}}
\frac{2^{k}\theta_{k}^{(s)}}{%
{\displaystyle\prod\limits_{j=0}^{k}}
\left(  2\theta_{j}^{(s)}+1\right)  }\right)  ,
\]
or using (\ref{theta})%
\begin{equation}
x_{n}^{(s)}=2^{-n}3^{\Theta\left(  n-1\right)  }\left(  s+%
{\displaystyle\sum\limits_{k=0}^{n-1}}
\frac{2^{k}\theta_{k}^{(s)}}{3^{\Theta\left(  k\right)  }}\right)  ,
\label{s1}%
\end{equation}
with%
\begin{equation}
\Theta\left(  x\right)  =%
{\displaystyle\sum\limits_{j=0}^{x}}
\theta_{j}^{(s)}. \label{Theta}%
\end{equation}
Setting $n=m$ and solving for $s$ in (\ref{s1}), we obtain $\left(
x_{m}^{(s)}=1\right)  $%
\begin{equation}
s=\frac{2^{m}}{3^{\Theta\left(  m-1\right)  }}-%
{\displaystyle\sum\limits_{k=0}^{m-1}}
\frac{2^{k}\theta_{k}^{(s)}}{3^{\Theta\left(  k\right)  }}. \label{s2}%
\end{equation}

Let $l=\Theta\left(  m-1\right)  .$ From (\ref{theta}) and (\ref{Theta}), we
see that $\Theta\left(  x\right)  $ is a step function with unit jumps at
$x=b_{1},b_{2},\ldots,b_{l},m,$ where $0\leq b_{1}<b_{2}<\cdots<b_{l}<m.$
Therefore, we can rewrite (\ref{s2}) as%
\[
s=\frac{2^{m}}{3^{l}}-%
{\displaystyle\sum\limits_{k=1}^{l}}
\frac{2^{b_{k}}}{3^{k}}.
\]
Finally, since $x_{m-3}^{(s)}=8,$ $x_{m-2}^{(s)}=4,$ $x_{m-1}^{(s)}=2$ and
$x_{m}^{(s)}=1,$ the penultimate jump must occur before or at $x=m-4.$ Thus,
$b_{l}\leq m-4$ and $l=\Theta\left(  m-1\right)  \leq m-3.$
\end{proof}

\begin{corollary}
\label{corollary}The Collatz conjecture is true if and only if every natural
number $n$ can be represented in the form%
\begin{equation}
n=\frac{2^{m}}{3^{l}}-%
{\displaystyle\sum\limits_{k=1}^{l}}
\frac{2^{b_{k}}}{3^{k}} \tag{C5}%
\end{equation}
for some $m,l,b_{1},\ldots,b_{l}\in%
\mathbb{N}
_{0},$ with
\[
0\leq l\leq m-3\ \ \text{and }\ 0\leq b_{1}<b_{2}<\cdots<b_{l}\leq m-4.
\]

\end{corollary}

Corollary \ref{corollary} is not a proof of the Collatz conjecture, but it
provides a lot of information on the set $\mathfrak{C}$ and the function
$\sigma_{\infty}(n).$ When $l=0,$ we recover the known fact that $2^{m}\in
S_{m},$ $\forall m\in%
\mathbb{N}
_{0}.$ For $l=1,$ we have the following result.

\begin{lemma}
For all $m\in%
\mathbb{N}
,$ we have%
\begin{align*}
\frac{2^{m}}{3}-\frac{2^{2k}}{3}  &  \in S_{m},\quad0\leq k\leq\frac{m-4}%
{2},\quad m\text{ \ even,}\\
\frac{2^{m}}{3}-\frac{2^{2k+1}}{3}  &  \in S_{m},\quad0\leq k\leq\frac{m-5}%
{2},\quad m\text{ \ odd.}%
\end{align*}

\end{lemma}

\begin{proof}
Let $n\in\Lambda_{m},$ with $l=1.$ We have
\[
n=\frac{2^{m}-2^{b_{1}}}{3}=2^{b_{1}}\times\frac{2^{m-b_{1}}-1}{3},\quad0\leq
b\leq m-4.
\]
Thus, $2^{m-b_{1}}\equiv1\left(  \operatorname{mod}2\right)  $ and therefore
$m-b_{1}\equiv0\left(  \operatorname{mod}2\right)  .$ Considering the cases
$m$ even and $m$ odd, the result follows.
\end{proof}

When $l=2,$ the situation is slightly more complicated. To simplify matters,
we restrict ourselves to the case of $n$ being odd.

\begin{proposition}
For all $m\geq5,$ with $m\neq6,8,$ we have%
\begin{align*}
\frac{2^{m}}{3^{2}}-\frac{2^{0}}{3}-\frac{2^{m-2-6k}}{3^{2}}  &  \in
S_{m},\quad1\leq k\leq\frac{m-2}{6},\quad m\geq10\text{ \ even,}\\
\frac{2^{m}}{3^{2}}-\frac{2^{0}}{3}-\frac{2^{m-4-6k}}{3^{2}}  &  \in
S_{m},\quad0\leq k\leq\frac{m-4}{6},\quad m\geq5\text{ \ odd.}%
\end{align*}

\end{proposition}

\begin{proof}
Let $n\in\Lambda_{m},$ odd, with $l=2.$ Then,
\[
n=\frac{2^{m}-3-2^{b_{2}}}{9}%
\]
and therefore%
\[
2^{m}-2^{b_{2}}=2^{b_{2}}\times\left(  2^{m-b_{2}}-1\right)  \equiv3\left(
\operatorname{mod}9\right)  .
\]
Considering all possible cases, we have

1) $2^{b_{2}}\equiv1\left(  \operatorname{mod}9\right)  $ and $2^{m-b_{2}%
}-1\equiv3\left(  \operatorname{mod}9\right)  $, which implies
\[
b_{2}\equiv0\left(  \operatorname{mod}6\right)  ,\quad m-b_{2}\equiv2\left(
\operatorname{mod}6\right)  .
\]

2) $2^{b_{2}}\equiv2\left(  \operatorname{mod}9\right)  $ and $2^{m-b_{2}%
}-1\equiv6\left(  \operatorname{mod}9\right)  $, which implies
\[
b_{2}\equiv1\left(  \operatorname{mod}6\right)  ,\quad m-b_{2}\equiv4\left(
\operatorname{mod}6\right)  .
\]

3) $2^{b_{2}}\equiv4\left(  \operatorname{mod}9\right)  $ and $2^{m-b_{2}%
}-1\equiv3\left(  \operatorname{mod}9\right)  $, which implies
\[
b_{2}\equiv2\left(  \operatorname{mod}6\right)  ,\quad m-b_{2}\equiv2\left(
\operatorname{mod}6\right)  .
\]

4) $2^{b_{2}}\equiv5\left(  \operatorname{mod}9\right)  $ and $2^{m-b_{2}%
}-1\equiv6\left(  \operatorname{mod}9\right)  $, which implies
\[
b_{2}\equiv5\left(  \operatorname{mod}6\right)  ,\quad m-b_{2}\equiv4\left(
\operatorname{mod}6\right)  .
\]

5) $2^{b_{2}}\equiv7\left(  \operatorname{mod}9\right)  $ and $2^{m-b_{2}%
}-1\equiv3\left(  \operatorname{mod}9\right)  $, which implies
\[
b_{2}\equiv4\left(  \operatorname{mod}6\right)  ,\quad m-b_{2}\equiv2\left(
\operatorname{mod}6\right)  .
\]

6) $2^{b_{2}}\equiv8\left(  \operatorname{mod}9\right)  $ and $2^{m-b_{2}%
}-1\equiv6\left(  \operatorname{mod}9\right)  $, which implies
\[
b_{2}\equiv3\left(  \operatorname{mod}6\right)  ,\quad m-b_{2}\equiv4\left(
\operatorname{mod}6\right)  .
\]
Thus, for $m$ even we shall have $m-b_{2}\equiv2\left(  \operatorname{mod}%
6\right)  $ or $b_{2}\equiv m-2\left(  \operatorname{mod}6\right)  $ and for
$m$ odd we need $m-b_{2}\equiv4\left(  \operatorname{mod}6\right)  $ or
$b_{2}\equiv m-4\left(  \operatorname{mod}6\right)  ,$ with $1\leq b_{2}\leq
m-4.$ Writing $b_{2}$ in terms of $m,$ the result follows.
\end{proof}

From Corollary \ref{corollary}, we can also get an idea of how the total
stopping time $\sigma_{\infty}(n)$ behaves if the Collatz conjecture is true.
Solving for $m$ in (C5) we have%
\[
m=\frac{1}{\ln\left(  2\right)  }\ln\left(  3^{l}n+3^{l}%
{\displaystyle\sum\limits_{k=1}^{l}}
\frac{2^{b_{k}}}{3^{k}}\right)  .
\]
In other words, $\sigma_{\infty}(n)$ lies on the family of parametric curves
\[
\frac{1}{\ln\left(  2\right)  }\ln\left(  3^{i}n+j\right)  ,\quad i,j\in%
\mathbb{N}
_{0},\quad i\leq j.
\]
For example, we have%
\[%
\begin{tabular}
[c]{|c|c|c|c|c|c|c|c|}\hline
$n$ & $2$ & $3$ & $4$ & $5$ & $6$ & $7$ & $8$\\\hline
$\sigma_{\infty}(n)$ & $1$ & $5$ & $2$ & $4$ & $6$ & $11$ & $3$\\\hline
$(i,j)$ & $\left(  0,0\right)  $ & $\left(  2,5\right)  $ & $\left(
0,0\right)  $ & $\left(  1,1\right)  $ & $\left(  2,10\right)  $ & $\left(
5,347\right)  $ & $\left(  0,0\right)  $\\\hline
\end{tabular}
.
\]

\begin{exercise}
Prove that
\[
\frac{\ln(n)}{\ln(2)}\leq\sigma_{\infty}(n)\quad\forall n\in%
\mathbb{N}
.
\]

\end{exercise}

\subsubsection{Binary sequences}

Another approach is to study the sequence $\left\{  \theta_{k}^{(s)}%
,\ k\geq0\right\}  ,$ which contains a wealth of information.

\begin{definition}
Let $\tau:\mathfrak{C}\rightarrow%
\mathbb{N}
$ be defined by%
\begin{equation}
\tau\left(  n\right)  =%
{\displaystyle\sum\limits_{k=0}^{\sigma_{\infty}(n)}}
\theta_{k}^{(n)}2^{k}. \label{tau}%
\end{equation}

\end{definition}

For example, we have%
\[%
\begin{tabular}
[c]{|c|c|c|c|c|c|c|c|c|}\hline
$n$ & $1$ & $2$ & $3$ & $4$ & $5$ & $6$ & $7$ & $8$\\\hline
$\tau(n)$ & $1$ & $2$ & $35$ & $4$ & $17$ & $70$ & $2199$ & $8$\\\hline
\end{tabular}
\ .
\]
Clearly, $\tau\left(  2^{n}\right)  =2^{n}$, $\forall n\in%
\mathbb{N}
_{0}.$

\begin{exercise}
Find $\tau(n)$ for $9\leq n\leq100.$
\end{exercise}

Let's study the image of $\Lambda_{m}$ by $\tau.$ We have%
\begin{align}
\tau\left(  \Lambda_{0}\right)   &  =\left\{  1\right\}  ,\quad\tau\left(
\Lambda_{1}\right)  =\left\{  2\right\}  ,\quad\tau\left(  \Lambda_{2}\right)
=\left\{  4\right\}  ,\quad\tau\left(  \Lambda_{3}\right)  =\left\{
8\right\}  ,\nonumber\\
\tau\left(  \Lambda_{4}\right)   &  =\left\{  16,17\right\}  ,\quad\tau\left(
\Lambda_{5}\right)  =\left\{  32,34,35\right\}  ,\quad\tau\left(  \Lambda
_{6}\right)  =\left\{  64,65,68,70\right\}  ,\label{taulambda}\\
\tau\left(  \Lambda_{7}\right)   &  =\left\{  128,130,136,137,140\right\}
,\quad\tau\left(  \Lambda_{8}\right)  =\left\{
256,257,260,272,274,280\right\}  .\nonumber
\end{align}

\begin{exercise}
Let
\[
\Lambda_{m}=\left\{  \lambda_{1}^{(m)},\ldots,\lambda_{N_{m}}^{(m)}\right\}
,
\]
where $N_{m}=\#\Lambda_{m}$ denotes the number of elements in the set
$\Lambda_{m}.$ Prove that $\forall m\in%
\mathbb{N}
_{0}$ there exist a sequence
\begin{equation}
2^{m}>\alpha_{1}^{(m)}\geq\cdots\geq\alpha_{N_{m}}^{(m)}=0, \label{alpha}%
\end{equation}
such that%
\[
\tau\left[  \lambda_{k}^{(m)}\right]  =2^{m}+\alpha_{k}^{(m)}.
\]

\end{exercise}

From (\ref{taulambda}), we have%
\[%
\begin{tabular}
[c]{|c|c|c|c|c|c|c|c|c|}\hline
$m$ & $1$ & $2$ & $3$ & $4$ & $5$ & $6$ & $7$ & $8$\\\hline
$\alpha_{1}^{(m)}$ & $0$ & $0$ & $0$ & $1$ & $3$ & $6$ & $12$ & $24$\\\hline
$\alpha_{2}^{(m)}$ & $0$ & $0$ & $0$ & $0$ & $2$ & $4$ & $9$ & $18$\\\hline
\end{tabular}
\ \ \ .
\]
It follows from (\ref{alpha}) that%
\[
\#\Lambda_{m}\leq\alpha_{1}^{(m)}+1,\quad\forall m\geq0.
\]

Using (\ref{s2}), we can define an inverse function for $\tau(n).$

\begin{definition}
Let $\phi:%
\mathbb{N}
\rightarrow%
\mathbb{Q}
$ be defined by%
\[
\phi(n)=\frac{2^{m}}{3^{\Phi\left(  m-1\right)  }}-%
{\displaystyle\sum\limits_{k=0}^{m-1}}
\frac{2^{k}\beta_{k}^{(n)}}{3^{\Phi\left(  k\right)  }},
\]
where $\beta_{m-1}^{(n)}\ldots\beta_{1}^{(n)}\beta_{0}^{(n)}$ is the binary
representation of $n,$ i.e.,
\[%
{\displaystyle\sum\limits_{k=0}^{m-1}}
2^{k}\beta_{k}^{(n)}=n\text{ \ and \ }\Phi\left(  x\right)  =%
{\displaystyle\sum\limits_{j=0}^{x}}
\beta_{j}^{(n)}.
\]
For example, we have%
\[%
\begin{tabular}
[c]{|c|c|c|c|c|c|c|c|c|}\hline
$n$ & $1$ & $2$ & $3$ & $4$ & $5$ & $6$ & $7$ & $8$\\\hline
$\phi(n)$ & $1$ & $2$ & $\frac{1}{3}$ & $4$ & $1$ & $\frac{2}{3}$ & $-\frac
{1}{9}$ & $8$\\\hline
\end{tabular}
\ \ .
\]

\end{definition}

\begin{exercise}
Find $\phi(n)$ for $9\leq n\leq100.$
\end{exercise}

It follows from Theorem \ref{Th} that
\[
\phi\circ\tau\left(  n\right)  =n,\quad\forall n\in%
\mathbb{N}
,
\]
while (\ref{alpha}) implies that
\[
S_{m}=\phi\left(  \left\{  2^{m},\ldots,2^{m}+\alpha_{0}^{(m)}\right\}
\right)  \cap%
\mathbb{N}
,\quad\forall m\geq0.
\]
In terms of $\phi,$ the Collatz conjecture reads%
\begin{equation}%
\mathbb{N}
\subset\phi\left(
\mathbb{N}
\right)  . \tag{C6}%
\end{equation}

With (C6), we finally reach a statement equivalent to the Collatz conjecture,
which is independent of the original formulation in terms of $T(x).$ Although
we have not succeeded in proving (C6), we hope that studying the function
$\phi(n)$ will shed new light on the Collatz problem.

\section{Further problems}

In the spirit of the Monthly, we offer a series of problems to the curious
reader. Those labeled "Exercise" are relatively easy to prove, "Problems"
denote results strongly supported by numerical evidence and "Conjectures" are
those that we would really wish to prove, but that may turn out to be false.

\begin{conjecture}
Prove that%
\begin{equation}
\sigma_{\infty}(n)\leq\delta(n)\ln(n)\quad\forall n\geq2, \label{C(n)}%
\end{equation}
where $\delta(n)$ is a slowly varying function, which might be eventually constant.
\end{conjecture}

\begin{definition}
The Abby-Normal numbers ($\mathcal{AN}$ numbers). Let the scaled total
stopping time $\gamma\left(  n\right)  $ be defined by
\[
\gamma\left(  n\right)  =\frac{\sigma_{\infty}(n)}{\ln(n)},\quad n\geq2.
\]
We say that $a_{k}$ is the $k-$th $\mathcal{AN}$ number if
\[
\gamma\left(  a_{k}\right)  =\max\left\{  \gamma\left(  n\right)
\ |\ a_{k-1}\leq n\leq a_{k}\right\}  ,\quad k\geq1
\]
with $a_{0}=2.$ In other words, $\left\{  \gamma\left(  a_{k}\right)
\right\}  $ is an increasing sequence of sharp lower bounds for the function
$\delta(n)$ defined in (\ref{C(n)}).
\end{definition}

\begin{exercise}
Show that%
\[%
\begin{tabular}
[c]{|c|c|c|c|c|c|c|c|}\hline
$n$ & $1$ & $2$ & $3$ & $4$ & $5$ & $6$ & $7$\\\hline
$a(n)$ & $3$ & $7$ & $9$ & $27$ & $230,631$ & $626,331$ & $837,799$\\\hline
$\gamma(n)\simeq$ & $4.55$ & $5.65$ & $5.92$ & $21.24$ & $22.51$ & $23.90$ &
$24.12$\\\hline
\end{tabular}
\ \ .
\]

\end{exercise}

From the results obtained by Eric Roosendaal \cite{Roosendaal}, it follows
that%
\begin{align*}
&  6,649,279,\quad8,400,511,\quad11,200,681,\quad15,733,191,\\
&  63,728,127,\quad3,743,559,068,799,\quad100,759,293,214,567,
\end{align*}
are possible $\mathcal{AN}$ numbers$.$ We have
\[
\gamma(100,759,293,214,567)\simeq35.17.
\]

\begin{exercise}
Find all $\mathcal{AN}$ numbers in the interval $[1,000,000,\quad
100,759,293,214,567].$
\end{exercise}

\begin{conjecture}
Prove that there exist infinitely many $\mathcal{AN}$ numbers.
\end{conjecture}

\begin{problem}
Let $V_{m}$ be the $m-$vector%
\[
V_{m}=\left[  x_{0}^{(m)},\ldots,x_{\sigma_{\infty}(m)}^{\left(  m\right)
}\right]  ,
\]
and $L:%
\mathbb{R}
^{N}\rightarrow%
\mathbb{R}
^{N}$ the linear operator defined by%
\[
L\left(  \left[  v_{1},v_{2},\ldots v_{N}\right]  \right)  =\left[
v_{2},v_{3},\ldots v_{N},v_{1}\right]  .
\]
Let $\theta\left(  m\right)  $ be the angle between $V_{m}$ and $L\left(
V_{m}\right)  .$ Prove that

\begin{description}
\item[(i)]
\[
\frac{1}{2}<\cos\left[  \theta\left(  m\right)  \right]  <\frac{7}{8}%
,\quad\forall m\in%
\mathbb{N}
.
\]

\item[(ii)]
\[
\frac{3}{4}<\cos\left[  \theta\left(  m\right)  \right]  <\frac{7}{8}%
,\quad\forall m\equiv1(2).
\]

\item[(iii)]
\[
\underset{k\rightarrow\infty}{\lim}\cos\left[  \theta\left(  a_{k}\right)
\right]  =\frac{7}{8}.
\]

\end{description}

Hint: See \cite{Gluck}.

\begin{problem}
Prove that:

\begin{description}
\item[(i)] $\forall k\geq6$ $\exists$ $m_{k}\in S_{k}$ such that $m_{k}+1\in
S_{k}.$

\item[(ii)] $\forall k\geq7$ $\exists$ $m_{k}\in S_{k}$ such that $m_{k}+2\in
S_{k}.$

\item[(iii)] $\forall l\in%
\mathbb{N}
$ $\exists$ $K\in%
\mathbb{N}
$ such that $\forall k\geq K$ $\exists$ $m_{k}\in S_{k}$ such that $m_{k}+l\in
S_{k}$ $.$
\end{description}

Hint: See \cite{Garner}.
\end{problem}
\end{problem}

\begin{exercise}
Prove that%
\[
\underset{k\rightarrow\infty}{\lim}\frac{\#S_{k+1}}{\#S_{k}}=\frac{4}{3}.
\]

\end{exercise}

\begin{exercise}
Let
\[
\zeta\left(  m\right)  =\#\left\{  2\leq k\leq m\ |\ \sigma_{\infty
}(k)=\ \sigma_{\infty}(k-1)\right\}  .
\]
Show that%
\begin{align*}
\zeta\left(  m\right)   &  =0,\ 2\leq m\leq12,\quad\zeta\left(  m\right)
=1,\ 13\leq m\leq14,\\
\zeta\left(  m\right)   &  =2,\ 15\leq m\leq18,\quad\zeta\left(  m\right)
=3,\ 19\leq m\leq20.
\end{align*}

\end{exercise}

\begin{conjecture}
Prove that%
\[
\underset{m\rightarrow\infty}{\lim}\frac{\zeta\left(  m\right)  }{m}=\frac
{1}{2}.
\]

\end{conjecture}

\begin{problem}
Prove that
\begin{equation}
\alpha_{1}^{(m)}=\left\lfloor 3\times2^{m-5}\right\rfloor ,\quad m\geq0
\label{alpha0}%
\end{equation}
and%
\[
\alpha_{2}^{(m)}=\left\lfloor 605\times2^{m-13}\right\rfloor ,\quad m\geq0,
\]
where $\left\lfloor \cdot\right\rfloor $ denotes the greatest integer function.
\end{problem}

\begin{problem}
Let
\[
S_{m}=\left\{  s_{1}^{(m)},\ldots,s_{N_{m}}^{(m)}\right\}  ,\quad
N_{m}=\#S_{m},
\]
with
\[
2^{m}=s_{1}^{(m)}>s_{2}^{(m)}>\cdots>s_{N_{m}}^{(m)}.
\]
Prove that%
\begin{align*}
\phi\left(  2^{m}+2^{2k}\right)   &  =s_{k+2}^{(m)},\quad0\leq k\leq\frac
{m-4}{2},\quad m\geq4,\quad m\equiv0(\operatorname{mod}2)\\
\phi\left(  2^{m}+2^{2k+1}\right)   &  =s_{k+2}^{(m)},\quad0\leq k\leq
\frac{m-5}{2},\quad m\geq5,\quad m\equiv1(\operatorname{mod}2).
\end{align*}

\end{problem}

\bibliographystyle{abbrv}
\bibliography{AIM}

\end{document}